\documentclass[11pt,reqno]{amsart}
\usepackage{amsmath}
\usepackage[curve,matrix,arrow]{xy}

\title{Tangential Quantum Cohomology of Arbitrary Order}

\author[S.\ J.\ Colley]{Susan Jane Colley}
\address{Department of Mathematics,
Oberlin College, Oberlin, Ohio 44074, USA}
\email{sjcolley@math.oberlin.edu}

\author[G.\ Kennedy]{Gary Kennedy}
\address{Ohio State University at Mansfield, 1680 University Drive,
Mansfield, Ohio 44906, USA}
\email{kennedy@math.ohio-state.edu}


\newtheorem{thm}[equation]{Theorem}
\newtheorem{lem}[equation]{Lemma}
\newtheorem{prop}[equation]{Proposition}

\theoremstyle{definition}
\newtheorem{defn}[equation]{Definition}

\theoremstyle{remark}

\numberwithin{equation}{section}

\newcommand{\mpsi}{\overline{\psi}} 
\newcommand{\xmod}[2]{\overline{M}_{0,#1}(X,#2)}
                                          
\newcommand{\pl}{\mathbf{P}^1} 
\newcommand{\pp}{\mathbf{P}^2} 

\newcommand{\gwc}[2]{\left\langle{#1}\mid{#2}\right\rangle}
\newcommand{\gwf}[4]{\left\langle\left\langle{#1}\mid{#2}\mid{#3}\mid{#4}\right\rangle\right\rangle}

\DeclareMathOperator{\vir}{vir}
\DeclareMathOperator{\id}{id}

\begin{document}

\begin{abstract}
 J.~Kock has previously defined a tangency quantum product on formal power series with coefficients in the cohomology ring of any smooth projective variety, and thus a ring that generalizes the quantum cohomology ring.  We further generalize Kock's construction by defining a $d$th-order contact product and establishing its associativity.
\end{abstract}

\maketitle

%
%
%
\section*{Introduction} \label{intro}

\par
A large chunk of the intersection theory of the moduli stacks of $n$-marked, genus 0 stable maps to a specified target variety $X$ is encoded by the quantum cohomology ring of $X$.  In particular, one can use the associativity of the ring structure to calculate Gromov-Witten invariants and, thus, in favorable circumstances, to calculate the number of rational curves in $X$ passing through an appropriate number of points.  In \cite{KockTh} and \cite{KockTQC}, J.~Kock defined a generalization of the quantum cohomology ring, called tangency quantum cohomology.  Kock's construction encodes a larger chunk of the intersection theory of the stacks of genus 0 stable maps than does ordinary quantum cohomology.  More specifically, it encodes the gravitational descendants or, to be more accurate, a modification of gravitational descendants (``enumerative descendants") better suited to questions arising in enumerative geometry, such as the computation of characteristic numbers of families of rational curves.  (An early avatar appeared in \cite{EK2}, which made a similar construction in the case of $\pp$.) 

\par
Kock's tangency quantum cohomology, however, deals only with descendants that allow at most a single tautological class (i.e., the first Chern class of the cotangent line bundle) to be specified at each of the $n$ marks.  Thus its enumerative applications (e.g., as in \cite{GKP}) are to situations involving ordinary tangency only.  We have discovered that his ideas may be generalized to deal with higher powers of tautological classes.  In this paper we thus complete the program begun by Kock:  we define a higher-order tangential quantum cohomology ring which allows one to deal with enumerative descendants up to an arbitrary specified order $d$.  We intend to explore the consequences of this structure for the purposes of enumerative geometry.  (The sort of application we envision would extend the ideas of our earlier papers \cite{ColleyK} and \cite{EK1}.)  In the present paper, however, we confine ourselves to constructing the ring and proving its associativity.

\par
In {\S}\ref{msosm} we begin with a brief reminder of basic notions:  moduli stacks of stable maps, tautological (psi) classes, and modified tautological classes.  We state Kock's ``key formulas" for the way in which modified tautological classes restrict to boundary components of the moduli space.  In {\S}\ref{pushpull} we present two push-pull results (Propositions \ref{pp} and \ref{ppp}) which are essential to the arguments in the sequel.  In {\S}\ref{GWclasses} we define Gromov-Witten and descendant classes, and present a series of technical lemmas for working with such classes.  Section \ref{tcp} contains our basic definition of the $d$th-order contact product and the main result (Theorem \ref{ringthm}), which states that this contact product gives an associative ring structure to a set of formal power series with coefficients in the Novikov ring of $X$.  The proof of Theorem \ref{ringthm} is devoted almost entirely to establishing associativity.  Finally in {\S}\ref{Ktcp} we show that when $d=1$ our contact product agrees with Kock's tangency quantum product.

\par
We remark on three aspects of our construction. First, although putative applications in enumerative geometry would employ a particular basis for the cohomology ring of the target space, our definitions and proofs have been formulated in a basis-free manner (except in the last section, where we compare our construction to Kock's).  For us, this recasting of Kock's constructions was an essential step in understanding how to craft a suitable generalization.  Second, we use the notion of virtual fundamental class; those whose interest is homogeneous spaces may simply use the fundamental class and ignore the more elaborate machinery. Third, it was already noted in \cite{EK2} and \cite{KockTQC} that the most naive generalization of the quantum product does not yield an associative product:  \cite{EK2} rectifies this by using an auxiliary space of ``stable lifts,'' whereas \cite{KockTQC} employs a ``twisted cup product.'' We likewise twist in a certain sense: note the factor $\exp(2z\chi_1)$ in Definition \ref{dbullet}. To our surprise, however, the twisting factor is no more complicated than that which already occurs at first order.

\smallskip
\par
We thank the following colleagues: Dan Abramovich for his help with, among other things, seeing how to avoid brutally long formulas by use of a suitably compact notation, Steve Kleiman and Ravi Vakil for their general advice and encouragement, Charles Cadman and Andrew Kresch for their technical advice.

\section{Moduli stacks of stable maps} \label{msosm}

\par
Let $X$ be any nonsingular complex algebraic variety, let $S$ be an (ordered) indexing set, and
let $\Delta \in H_2 (X)$ be the class of an effective curve or the zero class. (Throughout the paper,
we employ homology groups with rational coefficients.)  We will use the moduli stack $\xmod{S}{\Delta}$ of stable maps of genus $0$, as expounded in \cite{FultonP}, Chapter 7 of \cite{CoxKatz}, and Chapter 24 of \cite{Clay}.  If $S$ denotes the set $\{ 1, \ldots , n\}$, then we will also use the notation $\xmod{n}{\Delta}$.  The homology group $H_{2e}(\xmod{S}{\Delta})$ (where $e$ is the ``expected dimension") contains a virtual fundamental class $[\xmod{S}{\Delta}]^{\vir}$.  If $X$ is a homogeneous variety, then $e$ will coincide with the actual dimension of the moduli stack. For the rudiments of this theory, see \cite{BehrendFantechi}, \cite{LiTian1}, or \cite{Siebert}.

\par
If $S$ is the disjoint union of two subsets $S_1$ and $S_2$, and if
$\Delta=\Delta_1 + \Delta_2$, the associated \emph{boundary divisor} $D(S_1, \Delta_1 \mid S_2, \Delta_2)$
is the image of the morphism
\[ 
\rho_D \colon \xmod{S_1 \cup \{g_1\}}{\Delta_1} \times_X \xmod{S_2 \cup 
\{g_2\}}{\Delta_2} \to \xmod{S}{\Delta}
\]
obtained by gluing the marks labeled by $g_1$ and $g_2$.  (If $\Delta_i$ is the zero class we require that the cardinality of $S_i$ be at least 3.)
If $\Delta_1$ and $\Delta_2$ are both positive, we say (following Kock) that 
the corresponding $D(S_1, \Delta_1 \mid S_2, \Delta_2)$ is a 
\emph{hard boundary divisor}; if $\Delta_1$ or $\Delta_2$ is 
zero, we say it is a \emph{soft boundary divisor}.
 
\par
The ring $H^*(\xmod{S}{\Delta})$
contains \emph{evaluation classes}  $e_s^*\alpha$ pulled back from $X$ via the evaluation morphisms $e_s$ associated to the elements $s \in S$. It also contains \emph{tautological classes} defined as follows.
Consider the diagram
\[
\xymatrix{
\xmod{S\cup\{t\}}{\Delta} \ar@<3pt>[r]^<<<<<{\pi_{t}} \ar[d]_{e_{t}} & \xmod{S}{\Delta} 
\ar@<3pt>[l]^>>>>>{\sigma_s} \\
X & }
\]
in which $\pi_{t}$ denotes the morphism that forgets the extra mark,
and $\sigma_s$ is the section corresponding
to a particular $s \in S$.  The image of $\sigma_s$ is a soft
boundary divisor: it is the closure of the locus of maps
whose source curve consists of two twigs, one of which contracts to a
point and carries only the marks labeled by $s$ and $t$.  Denote the class of
this divisor by $D_{s,t}$.
If $\omega_{\pi_{t}}$ denotes the relative dualizing sheaf of 
$\pi_{t}$, then the \emph{cotangent line bundle} associated to $s$
is $\mathbb{L}_s := 
\sigma_s^*\omega_{\pi_{t}}$.  Its \emph{tautological class} (or \emph{psi class}) is
\[ 
\psi_s := c_1(\mathbb{L}_s) \in H^2(\xmod{S}{\Delta}).
\]

\par
For purposes of enumerative geometry, the tautological classes have an awkward feature: they are incompatible with pullbacks via forgetful morphisms. Instead,  we have 
\[ 
\pi_{t}^*\psi_s = \psi_s - D_{s, t}
\]
(where the tautological class on the left resides on 
$\xmod{S}{\Delta}$, while the one on the right resides on 
$\xmod{S\cup\{t\}}{\Delta}$).  To sidestep this problem, we work instead with 
the \emph{modified tautological class} $\mpsi_s$ defined (for $\Delta\neq 0$)
by
\begin{equation} \label{mpsidef}
\mpsi_s := \hat{\pi}_s^* \psi_s,    
\end{equation}
where $\hat{\pi}_s \colon \xmod{S}{\Delta} \to \xmod{\{s\}}{\Delta}$
is the morphism that forgets all marks except the one labeled by $s$.
(There is no need of this class when $\Delta= 0$.)

\par
Another important class is the \emph{diagonal class} $\delta_{r,s} \in 
H^2(\xmod{S}{\Delta})$, defined as the pull-back from $\xmod{\{r,s\}}{\Delta}$ of the 
Cartier divisor $D_{r,s}$.  Alternatively, $\delta_{r,s}$ is the sum of all 
boundary divisors representing maps in which the two associated marks lie on a contracting tree of twigs.

\begin{prop} \label{diagprops}
The following properties hold for the diagonal classes:
\begin{itemize}
\item[(a)]  
  $\pi_{r *}\left( \delta_{r,s}\cap[\xmod{S}{\Delta}]^{\vir}\right) 
  = [\xmod{S\setminus\{r\}}{\Delta}]^{\vir}$.
\item[(b)]
  $\delta_{r,s} \cup \delta_{r,t} = \delta_{r,s} \cup \delta_{s,t}$.   
\item[(c)]
  $\delta_{r,s} \cup e_r^*\alpha = \delta_{r,s} \cup e_s^*\alpha$.
\item[(d)]
  $\delta_{r,s} \cup \mpsi_{r} = \delta_{r,s} \cup \mpsi_{s}$. 
\end{itemize}    
\end{prop}
\begin{proof}
See Lemmas 1.3.2 and 1.3.3 of \cite{GKP}.
\end{proof}

\par
Our Theorem \ref{ringthm} relies on Kock's ``key formulas,'' which describe how the modified tautological classes on $\xmod{S}{\Delta}$ restrict to a boundary 
divisor $D(S_1, \Delta_1 \mid S_2, \Delta_2)$.  Let 
\[ 
j_D \colon 
\xmod{S_1 \cup \{g_1\}}{\Delta_1} \times_X \xmod{S_2 \cup 
\{g_2\}}{\Delta_2} 
\hookrightarrow 
\xmod{S_1 \cup \{g_1\}}{\Delta_1} \times \xmod{S_2 \cup 
\{g_2\}}{\Delta_2}
\]
denote inclusion of the fiber product into the Cartesian product. 

\newpage

\begin{thm} {\ }
\begin{enumerate}
    \item \textbf{Hard boundary case.}  
    If $\Delta_1, \Delta_2 >0$, then
    \begin{equation} \label{resthardeq}
    \rho_D^*\mpsi_s = 
    \begin{cases}
    j_D^*(\mpsi_s + \delta_{s,g_1}) 
	&\text{if $s \in S_1$} \\
    j_D^*(\mpsi_s + \delta_{s,g_2}) 
	&\text{if $s \in S_2$}    
    \end{cases}.    
    \end{equation}
    
    \item \textbf{Soft boundary cases.}
    \begin{enumerate}
      \item If $\Delta_1 > 0$, $\Delta_2 = 0$, then
         \begin{equation} \label{restsoft1eq}
	 \rho_D^*\mpsi_s = 
         \begin{cases}
          j_D^*\mpsi_s &\text{if $s \in S_1$} \\
          j_D^*\mpsi_{g_1} &\text{if $s \in S_2$}    
         \end{cases}.    
         \end{equation}
      \item If $\Delta_1 = 0$, $\Delta_2 > 0$, then
         \begin{equation} \label{restsoft2eq}
         \rho_D^*\mpsi_s = 
         \begin{cases}
         j_D^*\mpsi_{g_2} &\text{if $s \in S_1$} \\
         j_D^*\mpsi_s &\text{if $s \in S_2$}    
         \end{cases}.    
         \end{equation} 
    \end{enumerate}
\end{enumerate}    
\end{thm}    
\begin{proof}
For the hard boundary case, see Lemma 1.3.7 of \cite{GKP} or Section 1.4 of \cite{KockTQC}. Both cases are handled in Kock's thesis \cite{KockTh}.
\end{proof}

\section{Push-pull formulas} \label{pushpull}

\par
We begin this section with two technical lemmas in intersection theory. The basic background reference is \cite{Fulton}. We work in the context of bivariant rational equivalence theory, as expounded in \cite{BT}.  In subsequent sections, however, we will pass to singular homology and cohomology via the (Grothendieck) natural transformation. The paper \cite{Vistoli} (in particular its Section 5) explains how bivariant notions can be extended to the context of stacks.
\begin{lem} \label{claim0}
Consider the following fiber diagram of Deligne-Mumford stacks.
\[
\xymatrix{
M' \ar[d]_G \ar[r]^{i'} & N' \ar[d]^F \\
M \ar[r]^{i} & N }
\]
Suppose that $i$ is a regular embedding, $\gamma' \in A^*(N')$, and  $\nu' \in A_*(N')$. Then
\[ i^!(\gamma' \cap \nu') = i'^*\gamma' \cap i^!\nu'.  \]
\end{lem}
\begin{proof}
The Gysin pullback $i^!$ amounts to composition with $F^*[i]$, the pullback of the canonical orientation of the embedding (see \cite{BT}).  Thus, employing the notations of bivariant theory,
\[ i^!(\gamma' \cap \nu') = F^*[i] \cdot \gamma' \cdot \nu'. \]
In view of the remarks on page 652 of \cite{Vistoli}, we are working in the context of a commutative bivariant theory in the sense of \cite{BT}.  Hence this last expression equals
\[ i'^*\gamma' \cdot F^*[i] \cdot \nu'  = i'^*\gamma' \cdot i^! \nu' =  i'^*\gamma' \cap i^! \nu'. \]
\end{proof}
\begin{lem} \label{claim1}
Consider the following fiber diagram of Deligne-Mumford stacks.
\[
\xymatrix{
M'' \ar[d]_q \ar[r]^{i''} & N'' \ar[d]^p \\
M' \ar[d]_g \ar[r]^{i'} & N' \ar[d]^f \\
M \ar[r]^{i} & N }
\]
Let $i$ be a regular embedding, $\beta' \in A^*(N')$, and  $\nu'' \in A_*(N'')$. Then
\[ i^!(p^*\beta' \cap \nu'') = i''^*p^*\beta' \cap i^!\nu'' = q^*i'^*\beta' \cap i^!\nu''.  \]
\end{lem}
\begin{proof}
The first equality follows from Lemma \ref{claim0} with $\gamma' = p^*\beta'$, $G = g \circ q$, and $F = f \circ p$.  The second equality is obvious from the diagram.
\end{proof}

\par
The boundary divisors $D(S_1, \Delta_1 \mid S_2, \Delta_2)$ of $\xmod{S}{\Delta}$ are made up of fiber products of simpler moduli spaces. Thus our proofs in {\S}{\S}\ref{GWclasses}--\ref{tcp} use push-pull formulas, comparable to those in Proposition 1.7 of \cite{Fulton} and Lemma 3.9 of \cite{Vistoli}. Note that, since we assume that $X$ is a nonsingular variety, $A_*(X)$ and $A^*(X)$ are naturally isomorphic; we implicitly use this identification in these formulas.

\begin{lem}\label{avatar}
Suppose that $X$ is a nonsingular variety. Suppose that $M_1$ and $M_2$ are Deligne-Mumford stacks and that $g_1 \colon M_1 \to X$ is a proper morphism. Let $\alpha$ be a class in $A^*(M_1)$. Let $\mu_1$ be a class in $A_*(M_1)$ and $\mu_2$ a class in $A_*(M_2)$. In the fiber square
\[
\xymatrix{
& M_1\times M_2 \ar[dl]^{\id_{M_1} \times g_2} \ar[dr]_{g_1\times \id_{M_2}} \\
M_1 \times X \ar[dr]^{g_1\times \id_X} & & X \times M_2 \ar[dl]_{\id_X \times g_2} \\
& X \times X }
\]
we have
\[ \begin{split}
(g_1\times \id_{M_2})_{*}
\big( (\id_{M_1} \times g_2)^{*}&(\alpha \times 1) \cap (\mu_1 \times \mu_2) \big) \\
&=
(\id_X \times g_2)^{*}(g_1\times \id_X)_{*}
\big( (\alpha \times 1) \cap (\mu_1 \times [X]) \big)
\cap ([X] \times \mu_2).
\end{split} \]
\end{lem}
\begin{proof}
Both sides are equal to $g_{1*}(\alpha \cap \mu_1) \times \mu_2$.
\end{proof}
Note the hypothesis $\mu=\delta^{!}(\mu_1 \times \mu_2)$ in the following Proposition: in our later application this will motivate our definition of the virtual class on a fiber product of moduli spaces.
\begin{prop} \label{pp}
With the same hypotheses as in Lemma \ref{avatar}, consider 
the following fiber square. (Note that $q_2$ is proper.)
\[\xymatrix{
& M_1\times_X M_2 \ar[dl]^{q_1} \ar[dr]_{q_2} \\
M_1 \ar[dr]^{g_1} & & M_2 \ar[dl]_{g_2} \\
& X }
\]
Let $\mu=\delta^{!}(\mu_1 \times \mu_2) \in A_{*}(M_1\times_X M_2)$, where $\delta$ is the diagonal inclusion of $X$ into $X \times X$. Then
\[
q_{2*}\left(q_1^*\alpha\cap\mu \right)
=g_2^*g_{1*}(\alpha\cap \mu_1)\cap \mu_2.
\]
\end{prop}
\begin{proof}
Consider the following cube of morphisms, in which every face is a fiber square.
\begin{equation} \label{cube}
\xymatrix{
& {M_1\times_X M_2} \ar[rr]^{q_2} \ar[dl]_{q_1} \ar'[d]^{q_1}[dd] & & {M_2} \ar[dd]^{g_2} \ar[dl]|{(g_2,\id_{M_2})}\\
{M_1\times M_2} \ar[rr]^(.7){g_1\times\id_{M_2}} \ar[dd]_{\id_{M_1}\times g_2}& & {X\times M_2} \ar[dd]|(.3){\id_X\times g_2}\\
& {M_1} \ar'[r]^{g_1}[rr] \ar[dl]_{(\id_{M_1},g_1)}& & X \ar[dl]^{\delta}\\
{M_1\times X} \ar[rr]_{g_1\times\id_X} & & {X\times X}}
\end{equation}
Apply the Gysin homomorphism for $\delta$ to both sides of the formula in Lemma~\ref{avatar}.
On the left side we have
\begin{align*}
\delta^!(g_1\times \id_{M_2})_{*}
\big( (\id_{M_1} \times g_2)^{*}(\alpha \times 1) \cap (\mu_1 \times \mu_2) \big) 
&= q_{2*}\delta^!\big((\id_{M_1} \times g_2)^{*} (\alpha \times 1) \cap (\mu_1 \times \mu_2) \big) \\
&= q_{2*} \big(q_1^*(\id_{M_1}, g_1)^{*} (\alpha \times 1) \cap \delta^!(\mu_1 \times \mu_2) \big), \\
\intertext{by commuting the Gysin homomorphism with pushforward and then applying Lemma~\ref{claim1} to the diagram formed by the left and bottom faces of (\ref{cube}),}
&= q_{2*} \big(q_1^*\alpha \cap \delta^!(\mu_1 \times \mu_2) \big) = q_{2*}(q_1^*\alpha \cap \mu),
\end{align*}
since $(\id_{M_1},g_1)$ is the inclusion of the graph of $g_1$ in $M_1 \times X$.

\par
Now we apply $\delta^!$ to the right side of the formula in  Lemma ~\ref{avatar}:
\begin{align*}
\delta^!\big((\id_X \times g_2)^{*}(g_1\times \id_X)_{*}&\big( (\alpha \times 1) \cap (\mu_1 \times [X]) \big)
   \cap ([X] \times \mu_2)\big) \\
&= g_2^*\delta^!(g_1\times \id_X)_* \big( (\alpha \times 1) \cap (\mu_1 \times [X])\big)
   \cap \delta^!([X] \times \mu_2)\big).    \\
&= g_2^*\delta^!(g_1\times \id_X)_* \big( (\alpha \cap \mu_1)  \times [X]\big)
   \cap \delta^!([X] \times \mu_2)\big),      
\end{align*}
by Lemma \ref{claim1} applied to the following diagram derived from the right face of (\ref{cube}):
\[
\xymatrix{
M_2 \ar[d]_{g_2} \ar[rr]^(.4){(g_2, \id_{M_2})} & & X\times M_2 \ar[d]^{\id_X \times g_2} \\
X\ar[d]_{\id_X} \ar[rr]^(.4){\delta} & & X\times X \ar[d]^{\id_{X\times X}} \\
X \ar[rr]^(.4){\delta} & & X\times X }
\]
(Note that, in this case, $\delta^! = \delta^*$.)  Finally, employing commutativity of $\delta^!$ with pushforward, this last expression is
\[ g_2^*g_{1*}\delta^!\big( (\alpha \cap \mu_1) \times [X] \big) \cap \delta^!([X] \times \mu_2)
      = g_2^*g_{1*}(\alpha \cap \mu_1) \cap \mu_2,
\]
as desired.
\end{proof}

\begin{prop} \label{ppp}
With the same assumptions as in Proposition \ref{pp}, let $\beta$ be a class in $A^*(M_2)$.  Then we have
\[
q_{2*}\left(q_1^*\alpha \cup q_2^*\beta \cap \mu \right)
= g_2^*\left( g_{1*}(\alpha\cap \mu_1)\right) \cup \beta \cap \mu_2.
\]
\end{prop}

\begin{proof}
Apply the projection formula to the map $q_2$ and Proposition \ref{pp}.    
\end{proof}

\section{Gromov-Witten and descendant classes} \label{GWclasses}

\par
We now return to the context of homology and cohomology.
Suppose that $\gamma_1,\dots,\gamma_n$ are elements of $H^*(X)$ and that
$\Delta\in H_2(X)$ is the class of an effective curve (or the zero class).
We will employ the following notation adapted from Graber and Pandharipande, 
which we learned from Abramovich \cite{Ab}:  define a
\emph{Gromov-Witten class} by
\begin{equation}\label{gwclass}
\gwc{\gamma_1\cdot\ldots\cdot \gamma_n}{*}_{\Delta}
=(e_{n+1})_*\left( e_1^*\gamma_1 \cup \dots  \cup e_n^*\gamma_n
\cap[\xmod{n+1}{\Delta}]^{\vir}\right).
\end{equation}
This is intended to mimic the usual notation for a Gromov-Witten invariant, 
but note the formula defines a homology class on $X$ rather than a number.
More generally, we make the following definition.

\begin{defn} \label{gwdef}
A \emph{generalized Gromov-Witten class} is     
\begin{equation*}
\gwc{\gamma_1\cdot\ldots\cdot \gamma_n}{\gamma_{n+1}}_{\Delta}
=(e_{n+1})_*\left( \tilde{\gamma}_1 \cup \dots  \cup \tilde{\gamma}_n
\cup \tilde{\gamma}_{n+1}
\cap[\xmod{n+1}{\Delta}]^{\vir}\right),
\end{equation*}
where we use the convention that, for each $s=1,2,\dots,n+1$, either
\begin{equation*}
\gamma_s = \tilde\gamma_s \in H^*(\xmod{n+1}{\Delta})
\end{equation*}
or
\begin{equation*}
\gamma_s \in H^*(X) \text{ and } \tilde\gamma_s=e_s^*\gamma_s.
\end{equation*}
\end{defn}

\noindent
If $\gamma_{n+1}$ is the identity element, then we will revert 
to the notation of (\ref{gwclass}).  The following two lemmas follow 
immediately from the definitions.
\begin{lem} \label{lemmax}
If $\gamma_{n+1}\in H^*(X)$, then
\begin{equation*}
\gwc{\gamma_1\cdot\ldots\cdot \gamma_n}{\gamma_{n+1}}_{\Delta}
=
\gwc{\gamma_1\cdot\ldots\cdot \gamma_n}{*}_{\Delta}\cup\gamma_{n+1}.
\end{equation*}
\end{lem} 
\begin{lem} \label{lemmay}
If $\gamma_1,\gamma_2,\gamma_3\in H^*(X)$, then
\begin{equation*}
\gwc{\gamma_1\cdot\gamma_2}{\gamma_3}_0
=
\gamma_1\cup\gamma_2\cup\gamma_3.
\end{equation*}
\end{lem}

\par
To be able to work with arbitrary ordered index sets, we define a 
variant of Definition \ref{gwdef}:
\begin{equation} \label{dclass}
\gwc{\prod_{s\in S}\gamma_s}{\gamma_o}^{S\cup\{o\}}_{\Delta}
=(e_{o})_*\left(
{\bigcup_{s\in S} \tilde{\gamma}_s}
\cup  \tilde{\gamma}_o
 \cap[\xmod{S\cup\{o\}}{\Delta}]^{\vir}\right),
\end{equation}
with a similar convention and optional superscript.
(We use ``$o$" to suggest ``output.")
In the case where each $\tilde{\gamma}_s$
(as well as $\tilde{\gamma}_o$)
is of the form $(\mpsi_s)^t  \cup e_s^*\chi_s$
(where $t$ is a nonnegative exponent), we call the class 
in (\ref{dclass}) a \emph{descendant class}.

\par
Suppose that $\chi_0,\chi_1,\dots,\chi_d$ are classes in $H^*(X)$.
Let 
\begin{equation*}
\chi(s)=\sum_{t=0}^{d}(\mpsi_s)^t \cup e_s^* \chi_t.
\end{equation*}
(In the following section this will be used as a ``deforming element.'')
The following lemma is immediate from parts (c) and (d) of Proposition \ref{diagprops}.
\begin{lem}\label{swapchis}
$\delta_{r,s} \cup \chi(r) = \delta_{r,s} \cup \chi(s)$.
\end{lem}
\begin{lem}\label{vanish}
Suppose that $t$
is  a positive integer.
Then $\delta_{r,s}\cup(\mpsi_s)^t \cap [\xmod{S\cup\{r,s\}}{\Delta}]^{\vir} =0$.
\end{lem}
\begin{proof}
The indicated class is the pushforward of
\begin{equation*}
(\mpsi_s)^t \cap [\xmod{S\cup\{g_1\}}{\Delta}\times_{X}\xmod{\{g_2,r,s\}}{0}]^{\vir},
\end{equation*}
where $g_1$ and $g_2$ are the labels of the gluing marks.
The second factor is isomorphic to
\begin{equation*}
X\times\overline{M}_{0,\{g_2,r,s\}},
\end{equation*}
in which the second factor, the moduli space of stable rational curves 
with three marks, is simply a point. Hence $\mpsi_s=0$.
\end{proof}

\begin{lem} \label{replace}
Suppose that, for each $s\in T$, $\gamma_s$ is a class in $H^*(X)$. 
Suppose that $r$ is a nonnegative integer. Then
\begin{equation} \label{replaceright}
\gwc{\prod_{s\in T}\gamma_s\cdot\prod_{s\in U}\delta_{s,g}}{(\chi(g))^r}^{T\cup U\cup\{g\}}_{\Delta}
=\gwc{\prod_{s\in T}\gamma_s}{(\chi_0)^r}^{T\cup\{g\}}_{\Delta}.
\end{equation}
Similarly, for an arbitrary class $\zeta$,
\begin{equation} \label{replaceleft}
\gwc{(\chi(g))^r\zeta\cdot\prod_{s\in T}\gamma_s\cdot\prod_{s\in U}\delta_{s,g}}{\gamma_o}
^{\{g\}\cup T\cup U\cup\{o\}}_{\Delta}
=\gwc{(\chi_0))^r\zeta\cdot\prod_{s\in T}\gamma_s}{\gamma_o}
^{\{g\}\cup T\cup\{o\}}_{\Delta}.
\end{equation}
\end{lem}
\begin{proof}
We prove just (\ref{replaceright}); the proof of (\ref{replaceleft}) is similar.
The evaluation maps for the marks in $T$ (as well as the mark $g$) factor through
the map $f:\xmod{T\cup U\cup\{g\}}{\Delta} \to \xmod{T\cup\{g\}}{\Delta}$ which forgets the marks in $U$.
Furthermore, 
\begin{equation*}
f_*\left(\prod_{s\in U}\delta_{s,g}\cap[\xmod{T\cup U\cup\{g\}}{\Delta}]^{\vir}\right)=[\xmod{T\cup\{g\}}{\Delta}]^{\vir},
\end{equation*}
as one can see by induction, forgetting one mark at a time. (See formula (3), p.~168 of \cite{KockTQC} or Lemma~1.3.2 of \cite{GKP}.)  Thus, by the 
projection formula,
\begin{equation*}
\gwc{\prod_{s\in T}\gamma_s\cdot\prod_{s\in U}\delta_{s,g}}{(\chi(g))^r}^{T\cup U\cup\{g\}}_{\Delta}
=\gwc{\prod_{s\in T}\gamma_s}{(\chi(g))^r}^{T\cup\{g\}}_{\Delta}.
\end {equation*}
Now introduce an extra mark $h$. By the same argument (used in reverse),
we have
\begin{equation*}
\gwc{\prod_{s\in T}\gamma_s}{(\chi(g))^r}^{T\cup\{g\}}_{\Delta}
=\gwc{\prod_{s\in T}\gamma_s\cdot\delta_{h,g}}{(\chi(g))^r}^{T\cup\{h,g\}}_{\Delta}.
\end {equation*}
By Lemma \ref{swapchis}, the latter expression equals
\begin{equation*}
\gwc{\prod_{s\in T}\gamma_s\cdot\delta_{h,g}\,(\chi(h))^r}{*}^{T\cup\{h,g\}}_{\Delta}.
\end {equation*}
By Lemma \ref{vanish}, all contributions from positive powers of $\mpsi_h$ vanish.
Thus we may replace $\chi(h)$ by $\chi_0$, then again employ the swapping property of part (c) of Proposition \ref{diagprops} to obtain
\begin{equation*}
\gwc{\prod_{s\in T}\gamma_s\cdot\delta_{h,g}}{(\chi_0)^r}^{T\cup\{h,g\}}_{\Delta}.
\end {equation*}
Finally we obtain the right side of (\ref{replaceright}) by forgetting the mark $h$.
\end{proof}

\par
Our proof of associativity of the contact product in {\S}\ref{tcp} will be based on the following fiber diagram.
\begin{equation} \label{fiberdiag} 
\xymatrix@C-14pt{
 & \xmod{S_1\cup\{g_1\}}{\Delta_1}\times_{X} \xmod{S_2\cup\{g_2\}\cup\{o\}}{\Delta_2}
\ar[ddl] \ar[ddr]  \\
\\
  \xmod{S_1\cup\{g_1\}}{\Delta_1} \ar@/^/[dd]
\ar@/_/[dd]\ar[dd] \ar[ddr] & &
\xmod{S_2\cup\{g_2\}\cup\{o\}}{\Delta_2} \ar@/^/[dd] \ar@/_/[dd]
\ar[dd] \ar[ddl]  \\
\\
X & X & X  }
\end{equation}
Note that the evaluation map $e_s$ factors
through the first (respectively second) projection if $s\in S_1$ (respectively
$s\in S_2$). The evalution map $e_o$ likewise factors through the second projection.
There is also an evaluation map $e_g$ obtained by following around
either the left or right side of the diamond. (We use ``$g$" to suggest ``gluing.")

\par
The individual moduli spaces carry virtual fundamental classes based on perfect obstruction theories. Rather than grapple with the question of what is the appropriate obstruction theory on the fiber product, we simply declare that the virtual fundamental class of the fiber product is a Gysin pullback, via the diagonal inclusion of $X$ into $X \times X$:
\begin{align*}
[\xmod{S_1\cup\{g_1\}}{\Delta_1}\times_X \xmod{\{g_2\}\cup S_2\cup\{o\}}{\Delta_2}]^{\vir}&
\\
&\hspace{-1.3in}=\delta^!
\left(
[\xmod{S_1\cup\{g_1\}}{\Delta_1}]^{\vir}
\times
[\xmod{\{g_2\}\cup S_2\cup\{o\}}{\Delta_2}]^{\vir}
\right).
\end{align*}
Note that this is in accord with the hypothesis of Propositions \ref{pp} and \ref{ppp}.

\par
In this setup, we define
\begin{align*}
\gwf{\prod_{s\in S_1}\gamma_s}{*}{\prod_{s\in 
S_2}\gamma_s}{*}_{\Delta_1,\Delta_2}& \\
&\hspace{-1.3in}=(e_{o})_*\left(
{\bigcup_{s\in S_1 \cup S_2}\tilde {\gamma}_s}
 \cap[\xmod{S_1\cup\{g_1\}}{\Delta_1}\times_X \xmod{\{g_2\}\cup S_2\cup\{o\}}{\Delta_2}]^{\vir}\right)
\\
\intertext{and} 
\gwf{\prod_{s\in S_1}\gamma_s}{\gamma_g}{\prod_{s\in S_2}\gamma_s}{*}_{\Delta_1,\Delta_2}&\\
&\hspace{-1.3in}=(e_{o})_*\left(
{\bigcup_{s\in S_1 \cup S_2}\tilde {\gamma}_s}
\cup \tilde\gamma_g
 \cap[\xmod{S_1\cup\{g_1\}}{\Delta_1}\times_X \xmod{\{g_2\}\cup S_2\cup\{o\}}{\Delta_2}]^{\vir}\right),
\end{align*}
where our previous convention is now extended to three cases: 
either (1) $\tilde\gamma_s = e_s^*\gamma_s$, where $\gamma_s\in H^*(X)$, 
or (2) $\tilde\gamma_s$ is the pullback of a class $\gamma_s$ 
on the relevant factor of the fiber product, or (3) $\gamma_s = \tilde\gamma_s$, 
a class on the fiber product itself.

\par
\begin{lem}
In the first two cases (or any mixture of them),
 \begin{equation*}
\gwf{\prod_{s\in S_1}\gamma_s}{*}{\prod_{s\in S_2}\gamma_s}{*}_{\Delta_1,\Delta_2}
=\gwc{\gwc{\prod_{s\in S_1}\gamma_s}{*}_{\Delta_1}\cdot\prod_{s\in S_2}\gamma_s}{*}_{\Delta_2}.
\end{equation*}
Similarly in these first two cases
\begin{equation}\label{pupu}
\begin{split}
\gwf{\prod_{s\in S_1}\gamma_s}{\gamma_g}{\prod_{s\in S_2}\gamma_s}{\gamma_o}_{\Delta_1,\Delta_2}
&=\gwc{\gwc{\prod_{s\in S_1}\gamma_s}{\gamma_g}_{\Delta_1}\cdot\prod_{s\in S_2}\gamma_s}{\gamma_o}_{\Delta_2} \\
&=\gwc{\gwc{\prod_{s\in S_1}\gamma_s}{\gamma_g^{(1)}}_{\Delta_1}\gamma_g^{(2)}\cdot\prod_{s\in S_2}\gamma_s}{\gamma_o}_{\Delta_2},
\end{split}
\end{equation}
where $\gamma_g=\gamma_g^{(1)} \cup \gamma_g^{(2)}$ is any factoring.
\end{lem}

\begin{proof}
These follow from the push-pull Proposition \ref{ppp} together with repeated use of the projection formula.
\end{proof}

\section{The contact product} \label{tcp}

\par
As before, suppose that $\chi_0,\chi_1,\dots,\chi_d$ are classes in $H^*(X)$, and let 
\begin{equation*}
\chi(s)=\sum_{t=0}^{d}(\mpsi_s)^t \cup e_s^* \chi_t.
\end{equation*}
We call $\chi(s)$ a \emph{deforming element}.
Let $\alpha$ and $\beta$ be classes in $H^*(X)$. Working with moduli spaces 
having $n+3$ marks, we define the \emph{$d$th-order bullet product} 
(an auxiliary product that depends on the deforming element $\chi(s)$) by
\begin{equation} \label{dbullet}
\alpha\bullet\beta=\sum_{\substack{\Delta >0\\ n\ge 0}}\frac{z^n}{n!}q^{\Delta}
\gwc{\alpha\cdot\beta\cdot\prod_{s=1}^{n}\chi(s)}{\exp(2z\chi_1)}
^{\{a,b\} \cup \{1,2,\dots,n\} \cup \{c\}}_\Delta,
\end{equation}
where $q$ and $z$ are indeterminates and the sum is taken over all nonzero 
effective classes $\Delta \in H_2(X)$.   Note that $\alpha\bullet\beta$ is an element of 
$\Lambda[[z]]$, where $\Lambda$ denotes the Novikov ring of $X$ as defined in {\S}2.1 of \cite{Getzler}.  Extending this product by $\Lambda[[z]]$-linearity, 
we obtain a product on $\Lambda[[z]]$.  The \emph{$d$th-order contact 
product}, or \emph{$d$th-order tangential quantum product}, 
is then defined as
\begin{equation} \label{dcontact} 
\alpha * \beta =\alpha \cup \beta + \alpha \bullet \beta.
\end{equation}

\par
Our main result is as follows.
\begin{thm} \label{ringthm} 
The $d$th-order contact product is commutative and 
associative. The identity element $1 \in H^0(X)$ for the ordinary 
cup product is also the identity element for the contact product.
\end{thm}
We call the resulting ring 
the \emph{$d$th-order tangential quantum cohomology} and denote 
it by $Q^d H^*(X)$.

\par
\begin{proof}[Proof of Theorem \ref{ringthm}]
Commutativity is obvious. To verify that $1$ is the identity element, 
note that all classes appearing in
\begin{equation}\label{with1}
\gwc{1\cdot\beta\cdot\prod_{s=1}^{n}\chi(s)}{\exp(2z\chi_1)}
\end{equation}
are pullbacks via the morphism that forgets the first mark.
Since the fibers of this morphism have positive dimension, 
the projection formula implies that (\ref{with1}) vanishes. 
Thus $1\bullet\beta=0$ for all $\beta$.

\par
Let $\alpha, \beta, \gamma$ be classes in $H^*(X)$.
To verify associativity
we must show that
\begin{equation} \label{hopeqn}
(\alpha \bullet \beta )\bullet \gamma
    + (\alpha \bullet \beta )\cup \gamma
    + (\alpha \cup \beta )\bullet \gamma
=\alpha \bullet (\beta \bullet \gamma)
    + \alpha \cup (\beta \bullet \gamma)
    + \alpha \bullet (\beta \cup\gamma).
\end{equation}
To do so, we consider the moduli space
$\xmod{\{a,b,c,d\}\cup  \{1,\dots,n\}}{\Delta}$,
together with the forgetful morphism to the moduli space
$\overline{M}_{0,4}$ of stable genus $0$ curves with four marks.
On the latter space, which is isomorphic to $\pl$, let  $(a\, b \mid c\, d)$ denote
the point representing a two-twig curve carrying the points labeled $a$ and $b$ on one twig,
and those labeled $c$ and $d$ on the other. Let $\iota$ denote the inclusion of this point into
$\overline{M}_{0,4}$, and let $D(a\, b \mid c\, d)$ denote the fiber of the forgetful morphism
above this point. Define its virtual fundamental class via the Gysin pullback:
\begin{equation*}
[D(a\, b \mid c\, d)]^{\vir}:=
\iota^{!}[\xmod{\{a,b,c,d\}\cup  \{1,\dots,n\}}{\Delta}].
\end{equation*}
\par
Then in the homology group of $\xmod{\{a,b,c,d\}\cup \{1,\dots,n\}}{\Delta}$ we have the following equality:
\begin{equation*}
[D(a\, b \mid c\, d)]^{\vir}=
\sum_{\substack{\Delta_1,\Delta_2 \geq 0 \\ \Delta_1+\Delta_2 = \Delta \\ S_1\cup S_2=\{1,\dots,n\}}}
\rho_{*}[\xmod{\{a,b\}\cup  S_1 \cup \{g_1\}}{\Delta_1}
\times_X
\xmod{ \{g_2\}\cup \{c,d\}\cup S_2}{\Delta_2}]^{\vir},
\end{equation*}
where $\rho$ is the gluing map from the fiber product to $\xmod{\{a,b,c,d\}\cup  \{1,\dots,n\}}{\Delta}$.
This result, well-known when $X$ is a
convex variety, holds for smooth, projective $X$ from work of
K.~Behrend \cite{Behrend}, using techniques of
\cite{BehrendFantechi}.  (Also see \cite{Ab}.)
Since the points $(a\, b \mid c\, d)$ and $(b\, c \mid a\, d)$ are linearly equivalent, this virtual class
agrees with $[D(b\, c \mid a\, d)]^{\vir}$, from which we
deduce the following equation:
\begin{multline}  \label{holineqn}
\sum_{n\geq 0}\frac{z^n}{n!} \hspace{-.1in}
\sum_{\substack{\Delta_1,\Delta_2 \geq 0 \\
\Delta_1+\Delta_2 >0\\ S_1\cup S_2=\{1,\dots,n\}}}
\hspace{-.1in}q^{\Delta_1+\Delta_2}
\gwf{\alpha\cdot\beta\cdot\prod_{s\in S_1}\rho^*\chi(s)}{*}
{\gamma\cdot\prod_{s\in S_2}\rho^*\chi(s)}{\exp(2z\chi_1)}_{\Delta_1,\Delta_2} \\
 =
\sum_{n\geq 0}\frac{z^n}{n!}\hspace{-.1in}
\sum_{\substack{\Delta_1,\Delta_2 \geq 0 \\
\Delta_1+\Delta_2 >0\\ S_1\cup S_2=\{1,\dots,n\}}}
\hspace{-.1in}q^{\Delta_1+\Delta_2}
\gwf{\beta\cdot\gamma\cdot\prod_{s\in S_1}\rho^*\chi(s)}{*}
{\alpha\cdot\prod_{s\in S_2}\rho^*\chi(s)}{\exp(2z\chi_1)}_{\Delta_1,\Delta_2}.
\end{multline}
We will identify the left side of (\ref{holineqn}) as the left side of
(\ref{hopeqn}); that the right side of (\ref{holineqn}) agrees with the
right side of (\ref{hopeqn}) follows by permuting $a$, $b$, $c$, and
$d$ appropriately.

\par
In the hard boundary case, Kock's key formula (\ref{resthardeq}) tells us that if $s \in S_1$
(respectively $s \in S_2$),
then $\rho^*\chi(s)$ is the pullback from the first factor
(respectively second factor) of the class
\begin{equation*}
\chi(s)+\chi^\prime(s)\cup\delta_{s,g_1}
\quad\text{(respectively $\chi(s)+\chi^\prime(s)\cup\delta_{s,g_2}$),}
\end{equation*}
where
\begin{equation*}
\chi^\prime(s)=\frac{\chi(s)-e_s^*\chi_0}{\mpsi_s}
=\sum_{t=1}^{d}(\mpsi_s)^{t-1}\cup e_s^* \chi_t
\end{equation*}
and $\delta$ indicates the diagonal class of the two specified marks. 
(In the sequel we will omit the cup product symbols.)  Thus, by a binomial expansion, 
the hard boundary part of the left side of 
(\ref{holineqn})
is equal to
\begin{align*}
\sum_{n\geq 0} &\frac{z^n}{n!} \sum_{\substack{\Delta_1,\Delta_2 > 0 \\ 
T_1\cup U_1\cup T_2\cup U_2=\{1,\dots,n\}}}
q^{\Delta_1+\Delta_2} \\
& \gwf{\alpha\cdot\beta\cdot\prod_{s\in T_1}\chi(s)\cdot\prod_{s\in U_1}
\chi^\prime(s)\delta_{s,g_1}}{*}
{\gamma\cdot\prod_{s\in T_2}\chi(s)\cdot\prod_{s\in U_2}\chi^\prime(s)\delta_{s,g_2}}
{\exp(2z\chi_1)}_{\Delta_1,\Delta_2}. 
\end{align*}
We now rewrite this sum by citing many of our preliminary results. First,
since all of the classes involved are pullbacks from one factor or the other, 
by (\ref{pupu}) the sum equals
\begin{align*}
\sum_{n\geq 0} &\frac{z^n}{n!} 
\sum_{\substack{\Delta_1,\Delta_2 > 0 \\ T_1\cup U_1\cup T_2\cup U_2=\{1,\dots,n\}}}
q^{\Delta_1+\Delta_2} \\
&\gwc{\gwc
{\alpha\cdot\beta\cdot\prod_{s\in T_1}\chi(s)\cdot\prod_{s\in U_1}\chi^\prime(s)\delta_{s,g_1}}
{*}_{\Delta_1}
{{}\cdot\gamma\cdot\prod_{s\in T_2}\chi(s)\cdot\prod_{s\in U_2}\chi^\prime(s)\delta_{s,g_2}}}
{\exp(2z\chi_1)}_{\Delta_2}.
\end{align*}
Let $m_1$ indicate the cardinality of the set $U_1$, etc. Then by Lemma \ref{swapchis}
our sum equals
\begin{equation*}
\begin{split}
&\sum_{n\geq 0} \frac{z^n}{n!}
\sum_{\substack{\Delta_1,\Delta_2 > 0 \\ T_1\cup U_1\cup T_2\cup U_2=\{1,\dots,n\}}}
q^{\Delta_1+\Delta_2} \\
&\;\gwc{\gwc
{\alpha\cdot\beta\cdot\hspace{-.07in}\prod_{s\in T_1}\chi(s)\cdot\hspace{-.07in}\prod_{s\in U_1}\delta_{s,g_1}}
{(\chi^\prime(g_1))^{m_1}}_{\Delta_1}
{\hspace{-.15in}(\chi^\prime(g_2))^{m_2}\cdot\gamma\cdot\hspace{-.07in}\prod_{s\in T_2}
\chi(s)\cdot\hspace{-.07in}\prod_{s\in 
U_2}\delta_{s,g_2}}}{\exp(2z\chi_1)}_{\Delta_2}\hspace{-.05in}.
\end{split}
\end{equation*}
Using both parts of Lemma \ref{replace}, this is the same as
\begin{equation*}
\begin{split}
\sum_{n\geq 0} &\frac{z^n}{n!}
\sum_{\substack{\Delta_1,\Delta_2 > 0 \\ T_1\cup U_1\cup T_2\cup U_2=\{1,\dots,n\}}}
q^{\Delta_1+\Delta_2} \\
&\gwc{\gwc
{\alpha\cdot\beta\cdot\prod_{s\in T_1}\chi(s)}
{(\chi_1)^{m_1}}_{\Delta_1}
{(\chi_1)^{m_2}\cdot\gamma\cdot\prod_{s\in T_2}\chi(s)}}
{\exp(2z\chi_1)}_{\Delta_2}.
\end{split}
\end{equation*}
(Note the leading term of $\chi^\prime$ is $\chi_1$, and that we have replaced the class 
$\zeta$ of (\ref{replaceleft}) by a large expression.) Again by (\ref{pupu}), we obtain
\begin{align*}
\sum_{n\geq 0} &\frac{z^n}{n!} \sum_{\substack{\Delta_1,\Delta_2 > 0 \\
T_1\cup U_1\cup T_2\cup U_2=\{1,\dots,n\}}}
q^{\Delta_1+\Delta_2} \\
& \gwc{\gwc
{\alpha\cdot\beta\cdot\prod_{s\in T_1}\chi(s)}
{(\chi_1)^{m_1}(\chi_1)^{m_2}}
_{\Delta_1}
{{}\cdot\gamma\cdot\prod_{s\in T_2}\chi(s)}}
{\exp(2z\chi_1)}
_{\Delta_2}.
\end{align*}
Letting $n_1$ indicate the cardinality of $T_1$, etc., the last expression equals
\begin{equation*}
\begin{split}
&\begin{split}
\sum_{\substack{\Delta_1,\Delta_2 > 0 \\ n_1,n_2,m_1,m_2\geq 0}} 
  &\frac{z^{n_1+n_2}}{{n_1}!{n_2}!{m_1}!{m_2}!}q^{\Delta_1+\Delta_2} \\  
  &\gwc{\gwc{\alpha\cdot\beta\cdot\prod_{s\in T_1}\chi(s)}
    {(z\chi_1)^{m_1}(z\chi_1)^{m_2}}_{\Delta_1}
    {{}\!\!\cdot\gamma\cdot\prod_{s\in T_2}\chi(s)}}
    {\exp(2z\chi_1)}_{\Delta_2} 
\end{split}    \\
&\begin{split}
\quad=\sum_{\substack{\Delta_1,\Delta_2 > 0 \\ n_1,n_2\geq 0}} 
  &\frac{z^{n_1}}{{n_1}!}q^{\Delta_1}\frac{z^{n_2}}{{n_2}!}q^{\Delta_2} \\
  &\gwc{\gwc{\alpha\cdot\beta\cdot\prod_{s\in T_1}\chi(s)}
    {\exp(2z\chi_1)}_{\Delta_1}
   {{}\!\!\cdot\gamma\cdot\prod_{s\in T_2}\chi(s)}}
   {\exp(2z\chi_1)}_{\Delta_2} 
\end{split}   \\
&\quad=(\alpha\bullet\beta)\bullet\gamma.
\end{split}
\end{equation*}

\par
Now we deal with the soft boundary contributions on the 
left side of (\ref{holineqn}).  We will explain those terms in 
which $\Delta_2=0$; the terms in which $\Delta_1=0$ are handled 
by a completely analogous argument.  In this case, Kock's key formula 
(\ref{restsoft1eq}) implies that the deforming element restricts to a 
boundary divisor as follows:
\begin{equation*}
\rho^*\chi(s)=
\begin{cases}
\chi(s) &\text{if $s \in S_1$} \\
p^*\chi(g_1) &\text{if $s \in S_2$},  
\end{cases}
\end{equation*}
where $p$ indicates the projection of the fiber product to its first factor.
Thus the sum of the soft boundary contributions in which $\Delta_2=0$ is
\begin{align*} 
\sum_{n\geq 0} 
&\frac{z^n}{n!}\sum_{\substack{\Delta_1 > 0 \\  S_1\cup S_2=\{1,\dots,n\}}}
q^{\Delta_1}
\gwf{\alpha\cdot\beta\cdot\prod_{s\in S_1}\rho^*\chi(s)}{*}
{\gamma\cdot\prod_{s\in S_2}\rho^*\chi(s)}{\exp(2z\chi_1)}_{\Delta_1,0} \\
&=
\sum_{n\geq 0} \frac{z^n}{n!}
\sum_{\substack{\Delta_1 > 0 \\  S_1\cup S_2=\{1,\dots,n\}}}
q^{\Delta_1}
\gwf{\alpha\cdot\beta\cdot\prod_{s\in S_1}\chi(s)}{(\chi(g_1))^{|S_2|}}
{\gamma\cdot\prod_{s\in S_2}1}{\exp(2z\chi_1)}_{\Delta_1,0} \\
&=
\sum_{n\geq 0} \frac{z^n}{n!}
\sum_{\substack{\Delta_1 > 0 \\  S_1\cup S_2=\{1,\dots,n\}}}
q^{\Delta_1}
\gwc{\gwc{\alpha\cdot\beta\cdot\prod_{s\in S_1}\chi(s)}
{(\chi(g_1))^{|S_2|}}_{\Delta_1}
\cdot\gamma\cdot\prod_{s\in S_2}1}{\exp(2z\chi_1)}_0.
\end{align*}
All of the classes involved here are pullbacks via the forgetful morphism
forgetting $S_2$, but if this set is nonempty then its (virtual) 
fundamental class pushes forward to zero. Thus the only contributions 
to the sum are those in which $S_2$ is empty.  Hence it equals
\begin{equation*}
\sum_{n\geq 0} \frac{z^n}{n!}
\sum_{\Delta_1 > 0}
q^{\Delta_1}
\gwc{\gwc{\alpha\cdot\beta\cdot\prod_{s\in S_1}\chi(s)}
{*}_{\Delta_1}\cdot\gamma}{\exp(2z\chi_1)}
_0.
\end{equation*}
Successively using Lemma \ref{lemmax}, Lemma \ref{lemmay}, and again 
Lemma \ref{lemmax}, we see that this is the same as
\begin{align*}
\sum_{n\geq 0} \frac{z^n}{n!}
\sum_{\Delta_1 > 0}
q^{\Delta_1}&
\gwc{\gwc{\alpha\cdot\beta\cdot\prod_{s\in S_1}\chi(s)}
{*}_{\Delta_1}\cdot\gamma}{*}_0\cup
{\exp(2z\chi_1)}
 \\
&=
\sum_{n\geq 0} \frac{z^n}{n!} 
\sum_{\Delta_1 > 0}
q^{\Delta_1}
\gwc{\alpha\cdot\beta\cdot\prod_{s\in S_1}\chi(s)}
{*}_{\Delta_1}
\cup\gamma\cup
{\exp(2z\chi_1)}
 \\
&=
\sum_{n\geq 0} \frac{z^n}{n!}
\sum_{\Delta_1 > 0}
q^{\Delta_1}
\gwc{\alpha\cdot\beta\cdot\prod_{s\in S_1}\chi(s)}
{\exp(2z\chi_1)}_{\Delta_1}
\cup\gamma \\
&=(\alpha\bullet\beta)\cup\gamma.
\end{align*}
\end{proof}

\section{Kock's tangency quantum product} \label{Ktcp}

\par
In this concluding section, we show that, when $d=1$, our $d$th-order
contact product agrees with the tangency quantum product
defined by J.~Kock in \cite{KockTh} and \cite{KockTQC} and thereby
generalizes it.

\par
Given classes $\gamma_1,\ldots,\gamma_n \in H^*(X)$, T.~Graber, 
J.~Kock, and R.~Pandharipande (see \cite{GKP}) define (genus 0) 
\emph{enumerative descendants} 
$\left\{\bar{\tau}_{a_1}(\gamma_1)\cdots 
\bar{\tau}_{a_n}(\gamma_n)\right\}_{\Delta}$
as
\begin{equation} \label{endescdef} 
\left\{\bar{\tau}_{a_1}(\gamma_1)\cdots 
\bar{\tau}_{a_n}(\gamma_n)\right\}_{\Delta} :=
\int \mpsi_1^{a_1}\cup e_1^*\gamma_1 \cup \cdots\cup\mpsi_n^{a_n}\cup 
e_n^*\gamma_n \cap [\xmod{n}{\Delta}]^{\vir}.
\end{equation}
Our nonstandard use of curly braces, as opposed to angle brackets, is 
in order to distinguish the enumerative descendants (which are numbers) 
from the Gromov-Witten and descendant classes defined in {\S}\ref{GWclasses}.

\par
In \cite{KockTh} and \cite{KockTQC}, J.~Kock uses the enumerative 
descendants in (\ref{endescdef}) to define a type of 
``bullet product'' en route to defining his tangency quantum product.  
To do so, Kock begins with a basis $T_0,\ldots,T_r$ for $H^*(X)$.  
Let $x = x_0T_0+\cdots +x_rT_r$, $y = y_0T_0+\cdots +y_rT_r$ 
be arbitrary classes in $H^*(X)$, and define
\[ 
\gamma_{ij} := \int_X T_i \cup T_j \cup \exp(-2y) 
\qquad\text{and}\qquad 
(\gamma^{ef}) := (\gamma_{ij})^{-1}.
\]
(The $\gamma_{ij}$'s are coefficients for a deformation of the Poincar\'e 
metric.)  Define the generating function for the enumerative 
descendants by
\begin{align*}
\Gamma(x,y) &:= 
    \sum_{\Delta >0} q^{\Delta}\left\{\exp(\bar{\tau}_0(x) + 
	  \bar{\tau}_1(y))\right\}_{\Delta} \\
   &= \sum_{\substack{\Delta >0 \\ n \ge 0}} 
      \frac{q^{\Delta}}{n!} \int \bigcup_{s=1}^n 
      \left(e_s^*x + \mpsi_s \cup e_s^*y\right) 
      \cap [\xmod{n}{\Delta}]^{\vir}.
\end{align*}
Let $\Gamma_{ijl}$ denote the third-order partial derivative with 
respect to the $x$-variables:
\begin{equation} \label{Gammaderiv} 
\Gamma_{ijl} := \frac{\partial^3 \Gamma}{\partial x_i \partial x_j \partial x_l}
  = \sum_{\Delta > 0} q^{\Delta} 
  \left\{ \exp(\bar{\tau}_0(x) + \bar{\tau}_1(y))
  \bar{\tau}_0(T_i) \bar{\tau}_0(T_j) \bar{\tau}_0(T_l) \right\}_{\Delta}.
\end{equation}
(See {\S}3.1 of \cite{KockTQC}.)  Kock's bullet product 
$\bullet_{\text{Kock}}$ is defined with respect to the basis $T_0,\ldots,T_r$ by
\begin{equation} \label{KockBullet} 
   T_i \bullet_{\text{Kock}} T_j := 
   \sum_{l,m}\Gamma_{ijl}\gamma^{lm}T_m,
\end{equation}
which is an element of $\Lambda[[x_0,\ldots, x_r,y_0, \ldots,y_r]]$.
His tangency quantum product (see \cite{KockTQC}, {\S}3.5) is
\[
T_i *_{\text{Kock}} T_j := T_i \cup T_j + T_i \bullet_{\text{Kock}} T_j.
\]

\par
To compare Kock's bullet product with our own, which depends on a choice of deforming element
\[ 
\chi(s) = e_s^*\chi_0 + \mpsi_s \cup e_s^*\chi_1,
\]
we specialize Kock's indeterminates $x_0,\ldots, x_r, y_0, \ldots, y_r$ to a single indeterminate $z$ by requiring that
\begin{equation} \label{indeteqns}
\begin{aligned}
z\chi_0 = x_0 T_0 + \cdots + x_r T_r \\
z\chi_1 = y_0 T_0 + \cdots + y_r T_r.
\end{aligned} 
\end{equation}

\begin{prop} \label{KockequivProp}
    Under the specialization (\ref{indeteqns}), Kock's bullet product defined in (\ref{KockBullet}) agrees with 
    the $d=1$ case of the bullet product defined in (\ref{dbullet}).  
    Hence, when $d=1$, the $d$th-order contact product in (\ref{dcontact}) 
    is the same as Kock's tangency quantum product.
\end{prop}

\begin{proof}
Given any classes $\alpha, \beta  \in H^*(X)$, we show that 
$\alpha \bullet_{\text{Kock}} \beta = \alpha \bullet \beta$ by showing that 
\[
\int_X (\alpha \bullet_{\text{Kock}} \beta)\cup\exp{(-2y)}\cup \gamma 
= \int_X (\alpha \bullet \beta)\cup\exp{(-2y)}\cup \gamma,
\]
where $\gamma$ is an arbitrary class in $H^*(X)$.

\par
First, for all basis elements $T_i, T_j, T_k \in H^*(X)$, we have
\begin{alignat}{2}
\int_X (T_i &\bullet_{\text{Kock}}T_j) \cup \exp(-2y) \cup T_k
   & & \notag \\
&= \int_X\left(\sum_{l,m}\Gamma_{ijl}\gamma^{lm}T_m \right) \cup \exp(-2y) \cup T_k
   & & \quad\text{by (\ref{KockBullet}),} \notag \\
&= \sum_{l,m}\Gamma_{ijl}\gamma^{lm} \int_X \exp(-2y) \cup T_m \cup T_k  
   & & \notag \\
&= \sum_{l,m}\Gamma_{ijl}\gamma^{lm}\gamma_{mk} 
   & & \quad\text{by definition of $\gamma_{mk}$,} \notag \\
&= \Gamma_{ijk}
   & &\quad\text{since $(\gamma^{lm})$ and $(\gamma_{mk})$ are inverses.} \notag 
\end{alignat}    
Extending this result to arbitrary classes $\alpha, \beta, \gamma \in 
H^*(X)$ and using (\ref{Gammaderiv}), we have that
\begin{equation} \label{Kdegreeresult}
    \int_X (\alpha \bullet_{\text{Kock}}\beta ) \cup \exp(-2y) \cup \gamma
    = \sum_{\Delta >0}q^{\Delta}
      \left\{\exp(\bar{\tau}_0(x) + \bar{\tau}_1(y))
      \bar{\tau}_0(\alpha)\bar{\tau}_0(\beta)\bar{\tau}_0(\gamma)\right\}_{\Delta}.
\end{equation}

\par
Now we show that (\ref{Kdegreeresult}) holds when 
$\bullet_{\text{Kock}}$ is replaced by the $d=1$ case of the bullet product 
given in the definition of the tangential quantum product.  Note that under the specializations given by (\ref{indeteqns}), we have
\[ 
z\chi(s) = e_s^*x + \mpsi_s \cup e_s^*y   \qquad\text{and}\qquad   z\chi_1 = y.
\]
Thus
\begin{align*}
\int_X &(\alpha \bullet \beta ) \cup \exp{(-2y)} \cup \gamma \\
&= \int_X \sum_{\substack{\Delta >0 \\ n\ge 0}}\frac{z^n}{n!} q^{\Delta} 
  \gwc{\alpha \cdot\beta \cdot \prod_{s=1}^n 
  \chi(s)}{\exp{(2y)}}_{\Delta} \cup \exp{(-2y)} \cup \gamma \\
&= \sum_{\substack{\Delta >0 \\ n\ge 0}}
\frac{q^{\Delta}}{n!} 
\int_X \gwc{\alpha \cdot\beta \cdot \prod_{s=1}^n 
\left(e_s^*x + \mpsi_s \cup e_s^*y\right)}{*}_{\Delta} \cup \gamma \\
&= \sum_{\substack{\Delta >0 \\ n\ge 0}}\frac{q^{\Delta}}{n!} 
  \int a^*\alpha \cup b^*\beta \cup c^*\gamma \cup 
  \bigcup_{s=1}^n \left(e_s^*x + \mpsi_s \cup e_s^*y\right)  
  \cap[\xmod{n+3}{\Delta}]^{\vir}    \\
&= \sum_{\Delta >0}q^{\Delta}
  \left\{\bar{\tau}_0(\alpha)\bar{\tau}_0(\beta)\bar{\tau}_0(\gamma)
    \exp(\bar{\tau}_0(x) + \bar{\tau}_1(y))\right\}_{\Delta}.    
\end{align*}    
\end{proof}

\bibliographystyle{nyjalpha}
\ifx\undefined\bysame
\newcommand{\bysame}{\leavevmode\hbox to3em{\hrulefill}\,}
\fi

\end{document}